\documentclass[11pt]{article}
\usepackage{amssymb}

\def\marginpar#1{}

\makeatletter

\ifnum\@ptsize=0 \advance\hoffset by -0.3cm \fi

\ifnum\@ptsize=2 \advance\hoffset by 0.5cm \fi

\newtheorem{satz}{Satz}

\newtheorem{corollary}[satz]{Corollary}
\newtheorem{definition}[satz]{Definition}

\newtheorem{lemma}[satz]{Lemma}

\newtheorem{remark}[satz]{Remark}

\newtheorem{theorem}[satz]{Theorem}
\newtheorem{thevarthm}[satz]{\varthmname}

\newenvironment{varthm*}[1]{\trivlist\item[]{\bf #1.}\it}
   {\endtrivlist}

\newcommand\beginproof[1]{%
   \trivlist\item[\hskip\labelsep{\em #1.}]}
\newcommand\proof{\beginproof{Proof}}
\newcommand\proofof[1]{\beginproof{Proof of #1}}
\def\endproof{\hspace*{\fill}\endproofsymbol\endtrivlist}
\def\endproofsymbol{\frame{\rule[0pt]{0pt}{8pt}\rule[0pt]{8pt}{0pt}}}

\renewcommand\@seccntformat[1]{\csname the#1\endcsname.\enspace}
\renewcommand\epsilon{\varepsilon}

\renewcommand\phi{\varphi}
\def\vare{\varepsilon}

\newcommand\rounddown[1]{\left\lfloor#1\right\rfloor}

\newcommand\be{\begin{eqnarray*}}
\newcommand\ee{\end{eqnarray*}}

\newcommand\lra{\longrightarrow}
\newcommand\bbP{\mathbb P}

\newcommand\calo{{\cal O}}

\newcommand\refeq[1]{(\ref{#1})}

\newcommand\newop[2]{\def#1{\mathop{\rm #2}\nolimits}}

\newop\supp{supp}
\newop\AB{AB}
\newop\SB{SB}
\newop\End{End}
\newop\mult{mult}
\newop\NS{NS}
\newop\Nef{Nef}
\newop\lcd{lcd}

\def\be{\begin{eqnarray*}}
\def\ee{\end{eqnarray*}}
\def\ben{\begin{eqnarray}}
\def\een{\end{eqnarray}}

\def\qedsymbol{\frame{\rule[0pt]{0pt}{8pt}\rule[0pt]{8pt}{0pt}}}
\def\qed{\nopagebreak\hspace*{\fill}\qedsymbol\par\addvspace{\bigskipamount}}

\begin{document}
\noindent \parindent=0mm \parskip=.2cm

   \title{\large\bf An effective and sharp lower bound on Seshadri constants on surfaces with Picard number $1$}

   \author{\normalsize Tomasz Szemberg
      }

   \date{}

   \maketitle

   \thispagestyle{empty}
\bgroup\renewcommand\thesatz{\arabic{satz}}
\begin{abstract}
   On an algebraic surface with Picard number $1$ we compute
   in terms of
   the generator of the ample ray
   a lower bound for Seshadri constant valid at {\rm every}
   point of the surface. We show that this bound cannot be
   improved in general.
\end{abstract}
\section*{Introduction}
\label{intro}

   Seshadri constants were introduced by Demailly \cite{Dem92}.
   They measure the local positivity of an ample line bundle at a
   point. Though they are defined locally, they depend on
   the global geometry of the underlying variety and vice versa.

\begin{definition}
   Let $X$ be a smooth projective variety and $L$ an ample line
   bundle on $X$. Then
   $$\vare(L,x):=\inf_{x\in C}\frac{L.C}{\mult_x C}$$
   where the infimum is taken over all curves $C\subset X$ passing
   through $x$
   is the Seshadri constant at the point $x\in X$ (it is enough to
   consider irreducible curves).
\end{definition}

   By the ampleness criterion of Seshadri $\vare(L,x)$ is a
   positive real number. If $L$ is very ample, then it is easy to
   see that $\vare(L,x)\geq 1$ for all points $x\in X$.

   Shortly after Seshadri constants became an object of an
   independent study, Ein and Lazarsfeld \cite{EL93} proved a remarkably
   theorem that in most points an ample line bundle on a surface is locally as
   positive as a very ample one.

\begin{theorem}[Ein-Lazarsfeld]
   Let $S$ be a smooth projective surface and $L$ an ample line
   bundle on $S$. Then
   $$\vare(L,x)\geq 1$$
   for all points $x\in S$ away of at most countably many.
\end{theorem}

   On the other hand, Miranda \cite{PAG} provided examples showing
   that for any $\vare>0$ there exists a surface $S$, a point
   $x_0\in S$ and an ample line bundle $L$ on $S$ such that
   $\vare(L,x_0)<\vare$. In these examples the surfaces change as
   $\vare$ gets smaller and smaller. Moreover their Picard numbers
   grow reciprocally to $\vare$.

\section{The problem and the result}\label{problem}
   Ein and Lazarsfeld raised a natural question if there exists a
   single surface $S$ and sequences $L_n$ of ample line bundles
   and $x_n$ of points on $S$ such that
   $$\vare(L_n,x_n)\lra 0.$$
   This is not known up to now and it is conjectured that this is
   not possible i.e. that on a given surface there should be a
   universal lower bound on Seshadri constants of all ample line
   bundles.

   The result of Ein and Lazarsfeld was slightly improved by
   Oguiso \cite{Ogu02}.
\begin{theorem}[Oguiso]
   Let $S$ be a smooth projective surface and let $L$ be an ample
   line bundle on $S$. Then for an arbitrary $\delta>0$ the set of
   points $x$ such that
   $$\vare(L,x)\leq 1-\delta$$
   is finite.
\end{theorem}

   If the Picard number of $S$ is $1$, then there is essentially
   only the ample generator $L$ one has to take care of. In
   particular it follows from the above corollary that there
   exists a lower bound (namely the minimum over all points) for
   $\vare(L,x)$ but Oguiso theorem says nothing about estimating
   such a bound effectively.

\begin{corollary}\label{lower_non_eff}
   Let $S$ be a surface with Picard number $1$ with an ample
   generator $L$. Then there exists a number $\vare_0$ such that
   $$\vare(L,x)\geq\vare_0$$
   for all points $x\in S$.
\end{corollary}

   In order to make an effective statement one could revoke
   instead
   the big theorem of Matsusaka whose effective version on surfaces was
   proved by Fernandez del Busto \cite{FdB96}}.
\begin{theorem}[Fernandez del Busto]
   Let $L$ be an ample line bundle on a smooth projective surface
   $S$ with $a=L^2$ and $b=(K_S+L)L$. Then the line bundle $mL$ is
   globally generated (in particular $\vare(mL,x)\geq 1$) provided
   $$m>\frac{(b+1)^2}{2a}-1.$$
\end{theorem}

   Applying this result on a surface with Picard number $1$ yields
   the following effective statement.
\begin{corollary}
   Let $S$ be a smooth projective surface with Picard number $1$
   with an ample generator $L$ and let $r$ be an integer such that
   $K_S=rL$. Then
   $$\vare(L,x)\geq \frac{2L^2}{1+(r+4)^2(L^2)^2+2(r+3)L^2}$$
   for every point $x\in S$.
\end{corollary}

   This bound is hopelessly worse than the one stated in Theorem
   \ref{bound}.  In fact Theorem \ref{bound} is sharp and proving
   this constitutes the core of the present note.

\begin{theorem}\label{bound}
   Let $S$ be a smooth projective surface with $\rho(S)=1$ and let
   $L$ be an ample line bundle on $S$. Then for any point $x\in S$
   \begin{itemize}
   \item[(S)]
   $\vare(L,x)\geq 1$ if $S$ is not of general type and
   \item[(G)]
   $\vare(L,x)\geq\frac{1}{1+\sqrt[4]{K_S^2}}$ if $S$ is of
   general type.
   \end{itemize}
   Moreover both bounds are sharp.
\end{theorem}
\begin{remark}\rm
   It seems worth to note that on surfaces with Picard number $1$
   one has actually also a substantial improvement of the
   Ein-Lazarsfeld bound. Namely one has
   $$\vare(L,x)\geq \rounddown{\sqrt{L^2}}$$
   for $x$ general. This was observed by Steffens \cite{Ste98}.
\end{remark}
\proofof{Theorem \ref{bound}, case (S)}
   For the proof we go first through the Enriques-Kodaira
   classification of surfaces. Taking into account the assumption
   $\rho(S)=1$, there are only few cases.

   If $\kappa(S)=-\infty$, then $S=\bbP^2$ and it is well known
   that $\vare(\calo(1),x)=1$ for any point $x\in\bbP^2$. This
   verifies in particular that the bound stated in this part of the
   Theorem is sharp.

   If $\kappa(S)=0$, then $S$ is either abelian or K3.
   In the first case $S$ is a homogeneous variety, so $\vare(L,x)$
   does not depend on $x$ and by Ein-Lazasfeld Theorem we have
   $\vare(L,x)\geq 1$. Actually, Nakamaye \cite{Nak96} showed that
   if $\vare(L,x)=1$ on an abelian variety, then the variety is a
   product of an elliptic curve and a lower dimensional abelian
   variety. Note also that for abelian surfaces with Picard number $1$
   the exact values of Seshadri constants are known
   \cite{Bau98}.

   If $S$ is a K3 surface without $(-2)$-curves and $L$ is an ample
   line bundle on $S$, then $L$ is globally generated \cite{SD74}.
   This means that the morphism defined by the linear system $|L|$
   is finite, hence $\vare(L,x)\geq 1$ for all points $x\in S$.
\qed

\section{Seshadri constants of the canonical bundle}
   What remains are surfaces of general type. To complete the proof
   of Theorem \ref{bound} we need some preparations.

   First of all if the degree of the canonical divisor is not too small,
   then Reider's theorem \cite {Rei88} applies. More exactly we have the
   following lemma.

\begin{lemma}
   Let $S$ be a surface of general type with $\rho(S)=1$ (i.e.
   $K_S$ is ample) and $K_S\geq 5$. Then the bicanonical system
   $|2K_S|$ is base point free.
\end{lemma}
\proof
   This is just Reider's theorem for $K_S$. Note, that all
   exceptional cases in the theorem are immediately excluded under
   our assumptions.
\endproof
   As a corollary we get that
   \be\label{high}
   \vare(K_S,x)\geq\frac12
   \ee
   in this situation.

   Now we turn to the case $K_S^2\leq 4$. Then either $K_S$ is a
   primitive generator of the ample half-line or there exists an
   ample line bundle $L$ on $S$ with $K_S=2L$. In the latter
   situation it must be $L^2=1$ and consequently $K_S^2=4$.
   However such numerical invariants contradict the Riemann-Roch
   theorem for $L$. So we can assume that $K_S$ is a primitive
   line bundle. We obtain the following classification, which seems
   to be of independent interest.

\begin{lemma}\label{low}
   Let $S$ be a surface of general type with $\rho(S)=1$ and such
   that $K_S$ is primitive. Suppose
   that there exists a point $x\in S$ such that
   $$\vare(K_S,x)<1,$$
   then $K_S^2=1$, $q(S)=0$ and $p_g(S)\leq 2$ or $K_S^2=2$
   and $\vare(K_S,x)=\frac23$.
\end{lemma}
\proof
   By assumption we have that $K_S$ is ample.
   From \cite{CamPet90} it follows that there exists an irreducible curve
   $C\subset S$ such that $C$ computes the constant at $x$ i.e.
   $$\vare(K_S,x)=\frac{K_S.C}{m}<1$$
   with $m=\mult_x C$.
   There exists a positive integer $p$ such that
   $C\in|pK_S|$. So the above inequality yields $pK_S^2<m$ which
   is equivalent to $pK_S^2+1\leq m$.

   On the other hand, by the genus formula we have
   $$p_a(C)=1+\frac{p(p+1)}{2}K_S^2.$$
   A point of multiplicity $m$ causes the geometric genus of a
   curve to drop by at least $m\choose{2}$, so that
   $$1+\frac{p(p+1)}{2}K_S^2-{{m}\choose{2}}\geq 0.$$
   This gives $p^2(K_S^2)^2\leq 2+ p^2K_S^2$, which is possible
   only if $K_S^2=1$ or $K_S^2=2$ and $p=1$.

   By \cite[Theorem 11]{Bom73} $K_S^2=1$ implies $q(S)=0$.
   The inequality $p_g(S)\leq 2$ follows from the Noether
   inequality.

   If $K_S^2=2$, then $C$ is a canonical curve and in the
   inequalities above we have equality, so that in particular
   $m=3$. Then the Seshadri quotient is
   $$\frac{K_S.C}{m}=\frac23>\frac12.$$
\endproof

   Finally we take a closer look to surfaces of general type
   with $K_S^2=1$. They split
   again in two classes.

\subsection{Surfaces with $K_S^2=1$ and $p_q\leq 1$}
   If $p_g=0$ or $1$, then by the Riemann-Roch we have at least a pencil of
   bicanonical divisors. It is easy to check that the base locus
   of $|2K_S|$ in both cases consists only of points. Moreover for
   any point $x\in S$ there is an irreducible curve $D_x\in|2K_S|$
   passing through x. Let $C$ be any other irreducible curve on
   $S$ passing through $x$. Then we have
   $$2K_S.C=D_x.C\geq\mult_xC.$$
   This shows that the Seshadri quotient of $C$ satisfies
   $$\frac{K_S.C}{\mult_x C}\geq\frac12.$$
   The curve $D$ itself has arithmetic genus $4$, so it can have
   at most a triple point at $x$. Hence
   $$\frac{K_S.D_x}{\mult_xD_x}\geq\frac23>\frac12$$
   and we are done in these cases.

   Note that our argument is rather rough, in particular we didn't
   care if surfaces with given invariants and Picard number $1$
   exist. We will address this question later showing the optimality
   of the bound stated in the case (G) of Theorem \ref{bound}.

\subsection{Surfaces with $K_S^2=1$ and $p_g=2$}\label{pgdwa}
   The last case is that of a smooth surface $S$ of general type with $K_S^2=1$,
   $p_g(S)=2$ and $\rho(S)=1$.

   This time we can argue basically as in the preceding case but with
   the canonical pencil this time. This pencil consists of
   irreducible and reduced curves of genus $2$ all of whom pass
   through a single base point $x_0$ and meet there transversally.
   Let $x\in S$ be fixed and let $D_x$ be a curve in the pencil
   through $x$.
   If $C$ is an irreducible curve not in the pencil passing
   through $x$, then we have
   $$K_S.C=D_x.C\geq\mult_xC,$$
   so that in this case the Seshadri quotient is actually at least
   $1$. Now, it is not possible that all curves in the pencil are
   smooth. This can be seen either computing the topological Euler
   characteristic of the surface or with the argument that with
   all fibers smooth, the pencil would be an isotrivial family
   contradicting the assumption that $S$ is of general type. On
   the other hand, since the members of $|K_S|$ are curves of
   genus $2$ they can carry singularities with multiplicity at
   most $2$. We see that there must exist a canonical curve
   $D\in|K_S|$ and a point $x\in D$ with $\mult_xD=2$. Then
   $$\vare(K_S,x)=\frac{K_S.D}{\mult_xD}=\frac12.$$

   Summing up \refeq{high}, Lemma \ref{low} and the above discussion
   we have the following
\begin{varthm*}{Upshot}
   If $S$ is a surface of general type with $\rho(S)=1$, then
   $$\vare(K_S,x)\geq\frac12$$
   for arbitrary point $x\in S$.
\end{varthm*}
   At the end of the proof of case (G) of Theorem \ref{bound} we
   give an example showing that the bound in the Upshot is sharp.

\section{Primitive line bundles on surfaces of general type}
   In order to conclude the proof of Theorem \ref{bound} we have to
   study now the situation of $S$ being a surface of general type
   with Picard number $1$, $L$ an ample generator and $r$ a positive
   integer such that $K_S=rL$.

   From the Upshot stated above we get immediately a naive bound
   $$\vare(L,x)=\frac1r\vare(K_S,x)\geq\frac{1}{2r}\geq\frac{1}{2\sqrt{K_S^2}}$$
   but in fact we can do slightly better.

\proofof{Theorem \ref{bound}, case (G)}
   Assume that $x\in S$ is a point with a relatively low Seshadri constant
   (otherwise there is nothing to prove).
   Then there exists a curve $C\in|pL|$ computing this Seshadri
   constant
   $$\vare(L,x)=\frac{L.C}{m}=\frac{pL^2}{m},$$
   with $m=\mult_x C$.

   We have $p_a(C)=1+\frac12p(p+r)L^2$ and this gives an upper bound
   on the multiplicity $m$:
   $$m(m-1)\leq 2+p(p+r)L^2$$
   which is equivalent to
   $$m\leq \frac{1+\sqrt{9+4p(p+r)L^2}}{2}.$$
   Thus we have the following bound
   $$\vare(L,x)\geq\frac{2pL^2}{1+\sqrt{9+4p(p+r)L^2}}.$$
   The function on the right is growing for admissible values of $p$
   and $L^2$. Setting $L^2=1$ and $p=1$ we obtain
   \be\label{miedzy}
   \vare(L,x)\geq\frac{2}{1+\sqrt{13+\sqrt{K_S^2}}}.
   \ee
   Since the case of $K_S^2\leq 4$ was already discussed in Lemma
   \ref{low}, which in particular implies the bound stated in part
   (G) of the Theorem, we can assume that $K_S^2\geq 5$. But then it
   is easy to check that the number on the right in \refeq{miedzy}
   is greater or equal to our bound $\frac{1}{1+\sqrt[4]{K_S^2}}$ and
   this ends the proof of the inequality.

   Now we show that the bound is sharp.
   To this end let
   $S$ be a general surface of degree $10$ in the weighted projective
   space $\bbP(1,1,2,5)$. By adjunction we have that $K_S^2=1$.
   Moreover, sections of $K_S$ correspond to polynomials of degree
   $1$ in the weighted polynomial ring on $4$ variables. Thus
   $p_g(S)=2$ (see also \cite{Ste}).

   Steenbrink \cite{Ste} checked that a general surface $S$ of
   degree $d\geq 2+a+b$ in $\bbP(1,1,a,b)$ with $a$ and $b$ coprime has
   Picard number $1$. His result applies in our case. The existence
   of a point $x$ with $\vare(K_S,x)=\frac12$ follows now from the
   discussion in section \ref{pgdwa}.

   This example also shows that the bound given in the Upshot is in
   fact optimal.
\endproof
\section{Final remarks and a challenge}
   Looking back at the examples of Miranda we observe that in their
   case the lower bound of $\frac{1}{1+\sqrt[4]{|K_S^2|}}$ holds.
   This somehow gives a concrete effective number which could serve
   as a lower bound on arbitrary surface verifying in effect the
   conjecture stated in section \ref{problem}. It could be too much
   to state it as a conjecture but at least we dare a little
   challenge.
\begin{varthm*}{Question}
   Does there exist a (minimal) polarized surface $(S,L)$ and a point $x\in S$
   such that
   $$\vare(L,x)<\frac{1}{2+\sqrt[4]{|K_S^2|}}\;?$$
\end{varthm*}
   The appearance of $2$ in the above formulation accounts for
   the existence of Enriques surfaces which carry an ample line
   bundle $L$ with
   $\vare(L,x)=\frac12$, see \cite{Sze01}.

\begin{varthm*}{Acknowledgements}\rm\hspace{.1cm}
   Most of this work has been while the author visited University Duisburg-Essen.
   It is a pleasure to
   thank H\'el\`ene Esnault and Eckart Viehweg for their hospitality.
\end{varthm*}


\end{document}